# DEGENERATE DIFFUSIONS ARISING FROM GENE DUPLICATION MODELS

By Rick Durrett[1,2] and Lea Popovic[1,3]

*Cornell University and Concordia University*

We consider two processes that have been used to study gene duplication, Watterson's [*Genetics* **105** (1983) 745–766] double recessive null model and Lynch and Force's [*Genetics* **154** (2000) 459–473] subfunctionalization model. Though the state spaces of these diffusions are two and six-dimensional, respectively, we show in each case that the diffusion stays close to a curve. Using ideas of Katzenberger [*Ann. Probab.* **19** (1991) 1587–1628] we show that one-dimensional projections converge to diffusion processes, and we obtain asymptotics for the time to loss of one gene copy. As a corollary we find that the probability of subfunctionalization decreases exponentially fast as the population size increases. This rigorously confirms a result Ward and Durrett [*Theor. Pop. Biol.* **66** (2004) 93–100] found by simulation that the likelihood of subfunctionalization for gene duplicates decays exponentially fast as the population size increases.

**1. Introduction.** Studies have shown that a surprisingly large number of duplicated genes are present in all sequenced genomes, revealing that there is frequent evolutionary conservation of genes that arise through local events that generate tandem duplications, larger-scale events that duplicate chromosomal regions or entire chromosomes, or genome-wide events that result in complete genome duplication (polyploidization). Analyses of the human genome by Li (1980) have revealed that at least 15% of human genes are indeed duplicates, with segmental duplications covering 5.2% of the genome; see Bailey (2002). For more see the survey articles by Prince and Pickett (2002) and Taylor and Raes (2004).

Received January 2008; revised February 2008.

[1]Supported in part by an NSF/NIGMS grant.

[2]Supported by a grant from the probability program at NSF.

[3]Supported by NSF funded Focused Research Group.

Gene duplications are traditionally considered to be a major evolutionary source of new protein functions. The conventional view, pioneered by Ohno (1970) holds that gene duplication produces two functionally redundant, paralogous genes and thereby frees one of them from selective constraints. This unconstrained paralogue is then free to accumulate neutral mutations that would have been deleterious to a unique gene. The most likely outcome of such neutral evolution is for one of the paralogues to fix a null mutation and become a pseudogene, or to be deleted from the genome. A less frequently expected outcome is the fixation of mutations that lead to a new function. However, mutations that lead to such novel gene functions are extremely rare; see Walsh (1995), so the classical model predicts that few duplicates should be retained in genomes over the long term.

In Force (1999), introduced a new explanation of the preservation of duplicates called the duplication–degeneration–complementation (DDC) model. To explain, suppose that a gene performs two functions under the control of two different regulatory proteins that bind to the DNA upstream from the gene.

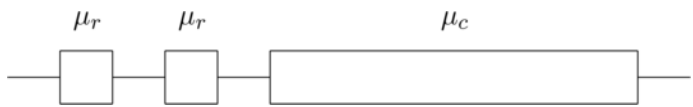

In the drawing the large rectangle is the gene, typically several hundred nucleotides, while the two small rectangles are the transcription factor binding sites, typically about 10 nucleotides. It is supposed that mutations which cause loss of a regulatory site happen at rate $\mu_r$ while those which cause the gene to completely lose function happen at rate $\mu_c$. In order to have the outcome called subfunctionalization in which the two genes specialize to do different functions, the first event must be a loss of a regulatory unit.

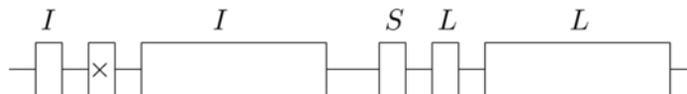

After this occurs, mutations in the indicated regions lead to inactivation of one gene copy, $I$, subfunctionalization, $S$, or are lethal $L$ since one of the functions is missing. It follows that the probability of subfunctionalization in a one lineage is

(1) $$\frac{4\mu_r}{4\mu_r + 2\mu_c} \cdot \frac{\mu_r}{2\mu_r + \mu_c} = 2\left(\frac{\mu_r}{2\mu_r + \mu_c}\right)^2.$$

If we make the simple assumption that $\mu_r = \mu_c$ this probability is 2/9, but if we observe that the gene region may easily be 30 times as large as the regulatory elements and set $\mu_c = 30\mu_r$ then the probability is 1/512.

Lynch and Force (2000) did simulations to assess the probability of subfunctionalization in a population with $N$ diploids. Taking $\mu_c = 10^{-5}$, which represents a gene with 1000 nucleotides and a per nucleotide per generation mutation rate of $10^{-8}$ they found that for $\mu_r/\mu_c = 3, 1, 0.3, 0.1$ the probability of subfunctionalization, $P(S)$ stayed constant at the value predicted by (1) until a population size of about 10,000 and then decreased rapidly to 0.

In a Research Experiences for Undergraduate project at Cornell, Ward and Durrett (2004) investigated the probability of subfunctionalization further using simulation and simple analytical computations. They found that if one plots the logarithm of $P(S)$ versus population size $N$, instead of Lynch and Force's plot of $P(S)$ versus $\log N$, one finds a straight line indicating exponential decay of $P(S)$. This conclusion holds for all four combinations of haploid and diploid organisms, and linked and unlinked duplicates.

Most multicellular organisms are diploid, but since simulations of Ward and Durrett (2004) show that the haploid and diploid situations are similar, we will stick with the simpler haploid case here. The linked case, where there is no recombination between the two gene copies, is relevant to tandem gene duplication in which unequal crossing over between two homologous chromosomes leads to a small region of the chromosome present in two copies. Ward and Durrett (2004) show that in the haploid linked case, the stochastic system converges to a system of differential equations, which always ends in a state in which one gene copy is lost. Since, for large $N$, the stochastic system is a small perturbation of this deterministic system, it follows from results of Friedlin and Wentzell (1988) that the probability of subfunctionalization decreases exponentially in $N$.

The unlinked case, in which there is free recombination between the two loci, occurs when copies reside on different chromosomes. This case is relevant to organisms such as the frog *Xenopus laevis* [Hughes and Hughes (1993)], teleost fishes [Bailey, Poulter and Stockwell (1978), Takahata and Maruyama (1979), Kimura (1980), Li (1980)] and yeast [Wolfe and Shields (1997), Seoighe and Wolfe (1998, 1999)], which have undergone whole genome duplication or tetraploidization, followed by a return to diploid inheritance. In these species, more than 25% of the duplicated genes have been preserved. Using simulation Ward and Durrett (2004) also made the startling discovery that in the haploid unlinked model the population frequencies of the various genetic states, which should be a six-dimensional diffusion, stayed near a one-dimensional curve of equilibria for the deterministic model.

The main point of this paper is to prove that observation. Our proof will show that in directions perpendicular to the curve of equilibria, the system is a small perturbation of a dynamical system with a strong push toward the curve. Since the point corresponding to subfunctionalization is not near this curve, an appeal to the results of Friedlin and Wentzell (1988) shows that the probability of subfunctionalization decays exponentially rapidly as the population size grows. This rigorously confirms the results Ward and Durrett (2004) found by simulation. Our proof gives no insight into the size of the constant in the exponential decay. Simulations suggest a value of 0.0005 for the diploid unlinked case, so for a population size in the millions, the probability of subfunctionalization will be extremely small.

A similar conclusion can be made about models in which the two gene copies control more than 2 functions which are regulated by different transcription factors. In the case of 3 functions the diffusion model becomes 14-dimensional. The equations for frequencies of all functional factors will again produce a one-dimensional curve of equilibria. However, since it becomes unwieldy to compute this curve or show that the 13-dimensional matrices associated with the linearization have all negative eigenvalues, we shall not try to prove this case.

As is often the case in work of this type, one can be certain of the truth of the mathematical fact, but cannot conclusively demonstrate its relevance to biology, because the result is based on assumptions which may or may not hold. In our case, some suspicious assumptions are:

1. constant population size,
2. mutations have no selective advantage,
3. the duplicated gene is present in all individuals.

Taking these in the order given, two referees have pointed out that the duplications that led to many gene families in Drosophila may have occurred well before the divergence of insects and chordates, perhaps in conditions of smaller population size. Of course, such gene families must be shared by all of these species, so this explanation applies only to the small number of families shared by all of these species. A second objection relevant to tetraploidization events is that these events initially involved only a single individual, so subfunctionalization may have occurred when the population is small. A third objection relevant to Drosophila is that they undergo large seasonal fluctuations and have spatial structure, so mutations may arise when the population size is large and become fixed in a subpopulation when it becomes small. Thus subfunctionalization can occur when the population size is small but this complements our result that it is unlikely when the population size is large.

As for the second assumption, we have assumed, following the proposers of the subfunctionalization explanation, that the mutations involved are neutral. In some cases, subfunctionalization means that the two genes are expressed in different tissues, and further mutations can improve the performance of the proteins. For a concrete example, consider the genes CYPB1 and CYPB3 of the Black Swallowtail studied by Wen et al. (2006). These genes, which have adapted to metabolize furanocoumarin found in different host plants are unlikely to be ancient duplicates common to all insects and chordates, but show a large number of amino acid altering mutations indicative of the action of positive selection (see Figure 2 in the paper cited).

Having admitted that positive selection is possible, we should also add that results of Walsh (1995) mentioned earlier show that positive adaption explains very few instances of gene duplication. As one referee pointed out, the probability of subfunctionalization may indeed be small, but there are many genes and many attempts, so there may be many successes anyway. One could object to this explanation by noting that if the subfunctionalization probability was 5% then there would be 19 failed duplicates in the genome for each success, and this is not observed. However, this objection vanishes if the duplications occur in a single individual and only become fixed in the population if it provides some benefit, since in this scenario failed attempts do not accumulate in the genome.

In summary, our result does not sound the death knoll for the popular concept of subfunctionalization. However, it does indicate that in order for it to be at work in a population which currently has a large population size, then it must be aided by positive selection or demographic factors.

## 2. Gene duplication models.
In this section we give precise statements on how degenerate diffusions arise in two distinct gene duplication models. Sections 3 and 4 contain proofs of Theorems 1 and 2 for Watterson's model in 2 dimensions. In Section 5 we give the proof of Theorem 3 for Lynch and Force's six-dimensional subfunctionalization model.

2.1. *Watterson's model.* To prepare for the complexities of the six-dimensional model, it is useful to start with Watterson's (1983) double recessive null model, which can be thought of as a special case of the subfunctionalization model in which the gene has only one function. Watterson considers a Wright–Fisher model with unlinked loci and diploid individuals. Unlinked loci means that the two gene copies are free to mutate independently and dependence only comes through viability of the gametes. In generation $n$ we have $2N$ letters that are either $A$ (working copy of gene 1) or $a$ (nonfunctional copy) and $2N$ letters that are either $B$ (working gene copy of gene 2) or $b$ (nonfunctional copy). To build up generation $n+1$ we repeat the following procedure until we have $N$ successes:

- Pick with replacement two copies of gene 1 and two copies of gene 2.
- An $A$, that is picked may mutate to $a$ with probability $\mu$. Likewise a $B$, that is picked may mutate to $b$ with probability $\mu$. The reverse mutation, which corresponds to undoing a specific deleterious mutation is assumed to have probability 0.

- If the result of our choices after mutation is not *aabb* then this is a success, and the new individual is added to the collection.

Note that after forming the new individual we do not keep track of the copies that reside in a single individual. This standard practice is referred to as the assumption of a "random union of gametes."

Letting $x$ and $y$ be the frequencies of $a$ and $b$ alleles, Kimura and King (1979) derived the diffusion limit of this model to be

$$(2) \qquad \begin{aligned} L_1 f = {} & \frac{1}{2} x(1-x) \frac{\partial^2 f}{\partial x^2} + 2N(1-x)(\mu - x^2 y^2) \frac{\partial f}{\partial x} \\ & + \frac{1}{2} y(1-y) \frac{\partial^2 f}{\partial y^2} + 2N(1-y)(\mu - x^2 y^2) \frac{\partial f}{\partial y}. \end{aligned}$$

To explain the coefficients in (2) we note that if we have a one locus Wright–Fisher model in which mutation from $A$ to $a$ occurs at rate $\mu$ and allele $a$ has a selective advantage $s$ then when time is run at rate $2N$ the diffusion approximation is

$$\frac{1}{2} x(1-x) \frac{\partial^2 f}{\partial x^2} + 2N[\mu(1-x) - sx(1-x)] \frac{\partial f}{\partial x}.$$

In our case an $A$ allele has fitness 1 and $a$ allele has fitness 0 if it is paired with another $a$ and two $b$'s, so the selective disadvantage of an $a$ allele is $s = xy^2$.

This diffusion limit in (2) is unusual since it does not assume that $s$ and $\mu$ are of order $1/N$. Since all individuals have fitness 1 or 0, it is not sensible to make this assumption about $s$, but one can, as Watterson did, assume $4N\mu \to \theta$. By using arguments that were clever but not completely rigorous, Watterson (1983) concluded that the mean time until loss of $A$ or $B$ had mean approximately

$$N[\log(2N) - \psi(\theta/2)],$$

where $\psi$ is the digamma function. Here we will give a simple proof of a version of his result. In our approach, we will assume that $\mu$ is a constant. Apart from simplifying the mathematics, our motivation is that while multicellular organisms have different population sizes, most have a mutation rate per nucleotide about $10^{-8}$, so if we assume genes have 1000 nucleotides, $\mu = 10^{-5}$.

To state our results, we begin by observing that solutions of the ordinary differential equation (ODE)

$$(3) \qquad \begin{aligned} \frac{dx}{dt} &= (1-x)(\mu - x^2 y^2), \\ \frac{dy}{dt} &= (1-y)(\mu - x^2 y^2) \end{aligned}$$

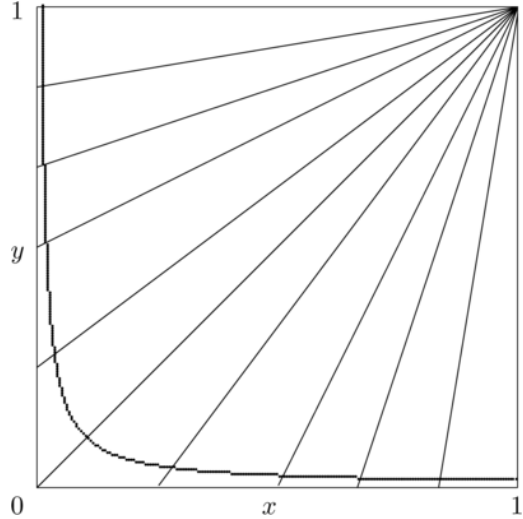

Fig. 1. *Picture of the state space for the Watterson model showing the curve of equilibria $x^2y^2 = \mu$ and the flow lines for the ordinary differential equation, which are the level curves for the projection map.*

have $(1-y_t)/(1-x_t)$ constant, so the solution moves along the line through $(1,1)$ until it hits the equilibrium curve $xy = \sqrt{\mu}$. Since the drift in the diffusion is $O(N)$ while its variance is only $O(1)$, it should not be surprising that the diffusion initially remains close to the ODE.

THEOREM 1. *Let $Z_t = (X_t, Y_t)$ be the diffusion in (2) and $Z_t^0$ the solution of (3) both starting from $(X_0, Y_0)$. Let $0 < \varepsilon < \sqrt{\mu}$. There is a constant $\gamma$, which depends on $\varepsilon$ so that if $N$ is large then for all $(X_0, Y_0) \in [\varepsilon, 1-\varepsilon]^2$*

$$E\left(\sup_{0 \leq t \leq \gamma \log N/N} |Z_t - Z_t^0|^2\right) \leq N^{-1/2}.$$

As the reader may have noticed this is not what the frequencies in the discrete Wright–Fisher model do on this time scale. They take significant jumps on the way to the curve. With more effort we could prove our results for the Wright–Fisher model, however, for simplicity, we will, as Watterson did, study the diffusion process.

With Theorem 1 established, it suffices to consider what happens when the diffusion starts near the curve. Given $(x,y)$ we define $(x^*, y^*)$ by $(1-y^*)/(1-x^*) = (1-y)/(1-x)$ and $x^*y^* = \sqrt{\mu}$. In words, $\Phi(x,y) = (x^*, y^*)$ is the point on the equilibrium curve which we would reach by flowing according to the ODE starting from $(x,y)$. See Figure 1 for a picture.

THEOREM 2. *Consider the diffusion $(X_t, Y_t)$ in (2). Let $\tau = \inf\{t : X_t = 1 \text{ or } Y_t = 1\}$ be the time to loss of $A$ or $B$. Let $0 < \delta < 1/2$. Suppose $|\mu - X_0^2 Y_0^2| \le N^{-\delta}$. Then if $N$ is large, with high probability we have $|\mu - X_t^2 Y_t^2| \le 2N^{-\delta}$ for all $t \le \tau$. Also, as $N \to \infty$ the process $X_t^* - Y_t^*$ converges in distribution to a diffusion process on the curve $\{xy = \sqrt{\mu}\}$, from which we obtain $E_0\tau \sim 2Nc_2(\mu)$.*

The coefficients of the limiting diffusion process can be explicitly calculated by applying Itô's formula to $h(\Phi(X_t, Y_t))$ where $h(x, y) = x - y$, but they have very complicated formulas. Figure 2 shows the drift $b(z)$, variance $a(z)$ of the limiting process and the quantity $-2b(z)/a(z)$ which appears in derivative of its natural scale

$$s'(y) = \exp\left(-\int_0^y -2b(z)/a(z)\,dz\right),$$

when $\mu = 10^{-4}$. The Green's function is shown in Figure 3. Recalling $E_0\tau = 2\int_0^{1-\mu} G(0, x)\,dx$ and integrating numerically we see that in our concrete case $c_2(\mu) = 6.993302$, so recalling how time is scaled we have $E_0\tau \approx 14N$ generations.

The result can be explained by noticing that the distance of the diffusion from the curve $\phi(X_t, Y_t) = (\mu - X_t^2 Y_t^2)^2$ is a Lyapunov function insuring exponential stability for the ODE system (3). Near the curve the stability of the drift beats the tendency of the diffusion to deviate from the curve, and forces $(X_t, Y_t)$ to follow where the flow would take it on the curve.

Watterson's argument proceeds by making a change of variables $\rho = 2N^{1/2}xy$ (we have $s = 1$ so his $S = N$) and

$$\eta = \begin{cases} \dfrac{x - y}{1 - y}, & \text{if } x \ge y, \\ \dfrac{x - y}{1 - x}, & \text{if } x \le y. \end{cases}$$

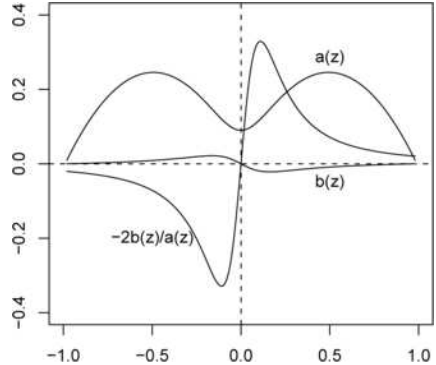

FIG. 2. *Diffusion coefficients for the limiting diffusion in the Watterson model.*

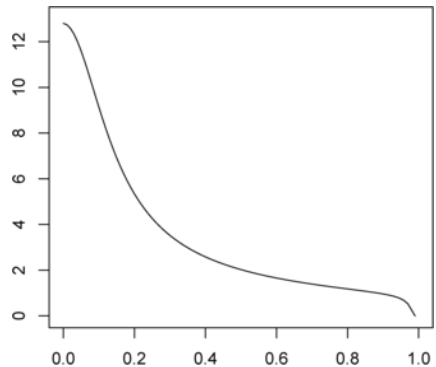

Fig. 3. *Green's function for the limiting diffusion in the Watterson model.*

Since his $\mu = \theta/4N$, $\rho^2 - \theta$ measure the displacement from the curve $x^2y^2 = \mu$, and $-1 \le \eta \le 1$ gives the position along the curve. The change of variables that defines $\eta$ is not twice differentiable when $x = y$ but by smoothing and passing to the limit this is not a problem.

The major unresolved difficulty with Watterson's argument is his assumption that $\rho$, which runs on a faster time scale than $\eta$, stays in equilibrium as $\eta$ evolves. A significant technical problem in justifying this is that when $\eta = 0$ the diffusion for $\rho$ does not have a stationary distribution. Indeed, it is the time spent with $\eta$ near 0 that results in the $N \log N$ asymptotics. Of course in his situation $\mu = O(1/N)$, so the situation is somewhat more complicated. Our proof, in contrast, uses the smooth change of variables $(x, y) \to x^* - y^*$, and to obtain our limit all we have to show is that the displacement from the curve stays close to 0.

2.2. *Subfunctionalization.* We turn now to the subfunctionalization model for unlinked loci in haploid organisms. In Watterson's model, the discrete-time Wright–Fisher model leads to a process that makes deterministic jumps off of the curve of equilibria. To avoid this, we will consider the Moran model version in which individuals are replaced one at a time, and in which mutation happens independent of reproduction. The biology behind the last choice is that we are considering the pool of gametes that exist at time $t$ in a population with overlapping generations.

To introduce our model, consider first an infinitely large population for which the allele frequency dynamics will be deterministic. Let $3 = 11$, $2 = 10$, $1 = 01$, and $0 = 00$ denote the four possible states of each gene copy, where 1 and 0 indicate presence or absence of the two functions, and let $x_i$ and $y_j$ denote the frequencies of states $i$ and $j$ at the first and second copy with $x_0 = 1 - x_3 - x_2 - x_1$, and $y_0 = 1 - y_3 - y_2 - y_1$. To simplify we will assume $\mu_r = \mu_c = b$. Let

(4)
$$w = x_3 + y_3 - x_3 y_3 + x_1 y_2 + x_2 y_1$$

be the mean fitness, that is, the probability the new individual chosen to replace the old one is viable. To explain the formula for $w$, we note that if either gene is in state 3, an event of probability $x_3 + y_3 - x_3 y_3$, then the offspring is always viable, whereas if neither gene is in state 3, the only viable combinations are $(1, 2)$ and $(2, 1)$. We are assuming the copies are unlinked so the events are independent.

The diffusion limit of this model is

$$L_2 f = \frac{1}{2} \sum_{i,j=1,2,3} x_i(\delta_{ij} - x_j) \frac{\partial^2 f}{\partial x_i \, \partial x_j} + \frac{1}{2} \sum_{i,j=1,2,3} y_i(\delta_{ij} - y_j) \frac{\partial^2 f}{\partial y_i \, \partial y_j}$$

(5)

$$+ \, 2NF(x_3, x_2, x_1, y_3, y_2, y_1) \cdot \nabla f,$$

where $F : \mathbb{R}^6 \mapsto \mathbb{R}^6$ is the vector field determining the evolution of the ODE system

$$dx_3/dt = -x_3 w + x_3 - 3bx_3,$$
$$dx_2/dt = -x_2 w + x_2(y_3 + y_1) + bx_3 - 2bx_2,$$
$$dx_1/dt = -x_1 w + x_1(y_3 + y_2) + bx_3 - 2bx_1,$$
(6)
$$dy_3/dt = -y_3 w + y_3 - 3by_3,$$
$$dy_2/dt = -y_2 w + y_2(x_3 + x_1) + by_3 - 2by_2,$$
$$dy_1/dt = -y_1 w + y_1(x_3 + x_2) + by_3 - 2by_1.$$

If we let $\alpha = 1 - 3b$, then the equations for $x_3$ and $y_3$ become

(7)
$$\frac{dx_3}{dt} = x_3(\alpha - w),$$
$$\frac{dy_3}{dt} = y_3(\alpha - w),$$

so the first and fourth equations for an equilibrium reduce to the single equation $w = \alpha$. Thus, if things are nice we will have a one-dimensional curve of equilibria.

To find one set of solutions we can begin by investigating the case in which $x_2 = x_1 = x$ and $y_2 = y_1 = y$ which gives us four equations

$$dx_3/dt = -x_3 w + x_3 - 3bx_3,$$
$$dx/dt = -xw + x(y_3 + y) + bx_3 - 2bx,$$
$$dy_3/dt = -y_3 w + y_3 - 3by_3,$$
$$dy/dt = -yw + y(x_3 + x) + by_3 - 2by,$$

which after some algebra can be explicitly solved. See Figure 4 for a graph of the solutions for $x_3, y_3, x, y$ in the special case $b = 0.001$.

It turns out that, in addition to the curve of equilibria, there are no other equilibria except for the fixed points of the dynamical system near $x_2 = 1$, $y_1 = 1$ and $x_1 = 1$, $y_2 = 1$, which correspond to subfunctionalization.

Equation (7) implies that if $y_3(0)/x_3(0) = r$ then we will have $y_3(t)/x_3(t) = r$ for all $t$. Given $z = (x_3, x_2, x_1, y_3, y_2, y_1)$, let $z^* = \Phi(z)$ be the point on the curve of equilibria with $y_3^*/x_3^* = y_3/x_3$. Numerical results show that starting from any $z$ the ODE will converge to $\Phi(z)$, but we do not know how to prove this for the dynamical system, so we will only consider the situation when the process starts close to the curve.

THEOREM 3. *Suppose that the mutation rate $b \leq 0.01$. Let $\tau = \inf\{t : X_0(t) = 1 \text{ or } Y_0(t) = 1\}$ be the time to loss of gene 1 or gene 2. Suppose $|Z_0 - Z_0^*| \leq 1/N^{1/4}$. Then if $N$ is large, with high probability we have $|Z_t - Z_t^*| \leq 2/N^{1/4}$ for all $t \leq \tau$. Also, when run at rate $2N$ the process $X_3^*(t) - Y_3^*(t)$ converges in distribution to a diffusion process on the curve of stable equilibria of (6), from which we get that $E_0 \tau \sim 2N c_3(b)$.*

The inspiration for our proof comes from Katzenberger (1991), who found conditions for a sequence of processes to converge to a diffusion on a submanifold of the original state space. While his result gives us the plan, it does not help much with the details. The hard work is to obtain enough information about the curve of equilibria to be able to check that the linearization of the drift in the five-dimensional space perpendicular to the curve has all eigenvalues with negative real part. (It is at this point we need the assumption that $b \leq 0.01$.) Linear ODE's of this type have Lyapunov functions, so patching these together and using computations in Section 3 of Katzenberger (1991) gives a Lyapunov function near the curve.

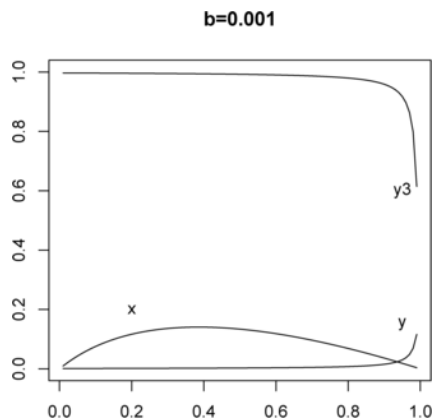

FIG. 4. *Picture of the solution of curve of equilibria in the subfunctionalization model where $x_2 = x_1 = x$, $y_2 = y_1 = y$ and the mutation rate has value $b = 0.001$ —the values of $(y_3, x, y)$ are shown against $x_3$ which is on the horizontal axis.*

Once we know the diffusion stays near the curve, the coefficients of the limiting diffusion can again be computed in terms of the derivatives of $\Phi(X_3, X_2, X_1, Y_3, Y_2, Y_1)$ and $h = x_3 - y_3$ by using Itô's formula. Figure 5 shows the drift $b(z)$, variance $a(z)$ of the limit of $X_3^* - Y_3^*$ and the quantity $-2b(z)/a(z)$ which appears in derivative of the natural scale when $b = 10^{-3}$. Even though the model is quite different, the curves are similar to those in Figure 2. The Green's function is shown in Figure 6. Integrating we see that in our concrete example $c_3(b) = 3.284906$, so for the process run at rate 1 we have $E_0 \tau \approx 6.5N$. In the case of *Drosophila* who have an effective population size of $N = 10^6$ this is 6.5 million generations or 650,000 years.

**3. Proof of Theorem 1.** We will prove this result by using a standard estimate from the Picard iteration proof of the existence of solutions of stochastic differential equations (SDE). We first consider a time rescaling $Z_t = (X_{t/2N}, Y_{t/2N})$, so that the coefficients of the SDE for $Z_t$ become

$$b(z) = ((1-x)(\mu - x^2 y^2), (1-y)(\mu - x^2 y^2)),$$
$$\sigma(z) = \begin{pmatrix} \sqrt{x(1-x)/2N} & 0 \\ 0 & \sqrt{y(1-y)/2N} \end{pmatrix}.$$

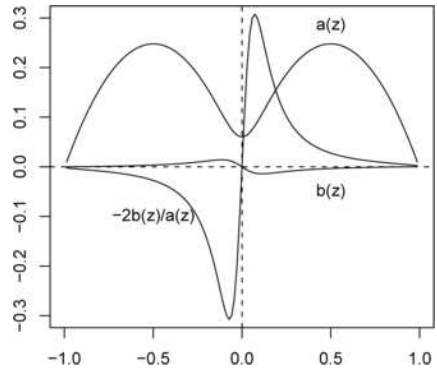

Fig. 5. *Diffusion coefficients for the limiting diffusion in the subfunctionalization model.*

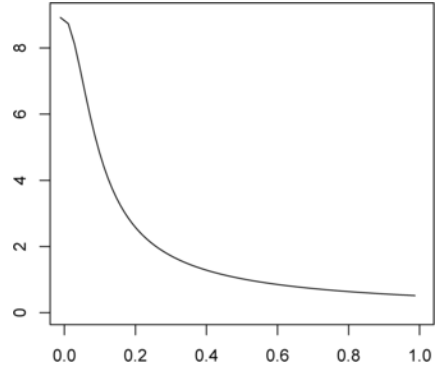

Fig. 6. *Green's function for the limiting diffusion in the subfunctionalization model.*

Let $b_\varepsilon(z)$ and $\sigma_\varepsilon(z)$ be the formulas that results when $x$ and $y$ are replaced by $(\varepsilon/2) \vee x \wedge (1 - \varepsilon/2)$ and $(\varepsilon/2) \vee y \wedge (1 - \varepsilon/2)$. We do this so that $b_\varepsilon$ and $\sqrt{2N}\sigma_\varepsilon$ are Lipschitz continuous with constant $K_\varepsilon$.

Let $Z_t^0$ be the solution of the ODE starting from $(X_0, Y_0) \in [\varepsilon, 1 - \varepsilon]^2$ with $\varepsilon < \mu$. Since along the boundary of the square, the drift points into the square except near the upper left and lower right corners, it is easy to see that

$$Z_t^0 - Z_0^0 = \int_0^t b_\varepsilon(Z_s^0) \, ds.$$

To begin the iteration to produce a solution of the modified SDE let

$$Z_t^1 = Z_0^0 + \int_0^t \sigma_\varepsilon(Z_s^0) \, dB_s + \int_0^t b_\varepsilon(Z_s^0) \, ds,$$

where $B_s$ is a standard two-dimensional Brownian motion. Since $u(1-u) \leq 1/4$ for all $u \in [0, 1]$, the $L^2$ maximal inequality for submartingales implies

$$E\left( \sup_{0 \leq s \leq t} |Z_s^1 - Z_s^0|^2 \right) \leq 4E|Z_t^1 - Z_t^0|^2 \leq t/N.$$

To continue we will use (2.3) on page 185 of Durrett (1996).

LEMMA 1. *If we let*

$$\tilde{U}_t = U_0 + \int_0^t \sigma_\varepsilon(U_s) \, dB_s + \int_0^t b_\varepsilon(U_s) \, ds,$$
$$\tilde{V}_t = V_0 + \int_0^t \sigma_\varepsilon(V_s) \, dB_s + \int_0^t b_\varepsilon(V_s) \, ds,$$

*then for* $B = (8T + 64)K_\varepsilon^2$

(8) $$E\left( \sup_{0 \leq s \leq T} |\tilde{U}_s - \tilde{V}_s|^2 \right) \leq 2E|U_0 - V_0|^2 + BE \int_0^T |U_s - V_s|^2.$$

Letting $Z_t^k = \tilde{Z}_t^{k-1}$, it follows by induction that for $k \geq 1$

$$E\left( \sup_{0 \leq s \leq T} |Z_s^k - Z_s^{k-1}|^2 \right) \leq \frac{B^{k-1}}{N} \frac{T^k}{k!}.$$

Using (2.6) on page 188 of Durrett (1996), we see that the solution $Z_t$ has

$$(9) \quad E\left(\sup_{0 \le s \le T} |Z_s - Z_s^0|^2\right) \le \left(\sum_{k=1}^{\infty} \sqrt{\frac{B^{k-1}}{N} \frac{T^k}{k!}}\right)^2 \le \frac{1}{N}\left(\sum_{k=1}^{\infty} \frac{((8T+64)K_\varepsilon)^k}{\sqrt{k!}}\right)^2.$$

The form of the right-hand side derives from the fact that it is the $L^2$ norm $\|X\|_2 = (EX^2)^{1/2}$ that satisfies the triangle inequality. The right-hand side increases rapidly with $T$, so we take $T$ to be a constant to end up with a bound of the form $C_T/N$.

To extend the comparison to longer intervals of time, let $Z_t^{1,0}$ and $Z_t^{1,1}$ be solutions of the ODE and SDE starting from $Z_T^0$, and let $Z_t^{1,2}$ be the solution of the SDE starting from $Z_T$. It follows from (9)

$$E\left(\sup_{0 \le s \le T} |Z_s^{1,1} - Z_s^{1,0}|^2\right) \le C_T/N.$$

Using (8) with Gronwall's inequality [see (2.7) on page 188 of Durrett (1996)]

$$E\left(\sup_{0 \le s \le T} |Z_s^{1,2} - Z_s^{1,1}|^2\right) \le \frac{2C_T}{N}e^{BT}.$$

Combining the two estimates using the triangle inequality gives

$$E\left(\sup_{0 \le s \le T} |Z_s^{1,2} - Z_s^{1,0}|^2\right) \le \frac{C_T}{N}(1 + \sqrt{2e^{BT}})^2.$$

To continue for $k \ge 2$ let $Z_t^{k,0}$ and $Z_t^{k,1}$ be solutions of the ODE and SDE starting from $Z_{kT}^0 = Z_T^{k-1,0}$, and let $Z_t^{k,2}$ be the solution of the SDE starting from $Z_T^{k-1,2}$. We will show by induction that

$$(10) \quad E\left(\sup_{0 \le s \le T} |Z_s^{k,2} - Z_s^{k,0}|^2\right) \le \frac{C_T}{N}\left(\sum_{j=0}^{k} (2e^{BT})^{j/2}\right)^2.$$

Our previous result shows that this hold for $k = 1$. It follows from (9)

$$E\left(\sup_{0 \le s \le T} |Z_s^{k,1} - Z_s^{k,0}|^2\right) \le C_T/N.$$

Using (8) with Gronwall's inequality, and our induction assumption we get

$$E\left(\sup_{0 \le s \le T} |Z_s^{k,2} - Z_s^{k,0}|^2\right) \le \frac{2C_T}{N}\left(\sum_{j=0}^{k-1} (2e^{BT})^{j/2}\right)^2 e^{BT}.$$

Combining the two estimates gives (10). If we choose $T$ so that $2e^{BT} = 4$, then

$$\left(\sum_{j=0}^{k} (2e^{BT})^{j/2}\right)^2 \le \frac{(2e^{BT})^k}{[1 - (2e^{BT})^{-1/2}]^2} = 4^{k+1}$$

and summing over $k = 1, \ldots, m$ gives

$$E\left(\sup_{0 \le s \le mT} |Z_s - Z_s^0|^2\right) \le \frac{C_T}{N} \frac{4^{m+1}}{1 - 1/4} = \frac{16}{3} C_T N^{-2/3} \le N^{-1/2},$$

when we choose $m = (\log N)/(3 \log 4)$. This proves Theorem 1 with $\gamma = T/(6 \log 4)$ and time rescaled back to its original rate.

**4. Proof of Theorem 2.** The proof of Theorem 2 is given in four steps. We first find the map $\Phi(x, y)$ that takes a point $(x, y)$ to where the ODE trajectory would take it on the curve; and apply Itô's formula to find the SDE of $\Phi(X_t, Y_t)$. We then use the exponential Lyapunov function for the ODE to show that the diffusion stays close to the process $\Phi(X_t, Y_t)$; and finally use its SDE to obtain the limiting process.

4.1. *Projection map of the trajectory.* Given $x, y$ in the unit square the ODE will take it to the point $x^*, y^*$ such that $x^* y^* = \sqrt{\mu}$ while $(1 - y)/(1 - x) = (1 - y^*)/(1 - x^*)$. Solving for $x^*$ gives a quadratic equation

$$(1 - y)(x^*)^2 - (x - y)x^* - \sqrt{\mu}(1 - x) = 0$$

and the root we want is

$$x^* = \frac{1}{2}\left(\frac{x - y}{1 - y} + \sqrt{\frac{(x - y)^2}{(1 - y)^2} + 4\sqrt{\mu}\frac{1 - x}{1 - y}}\right).$$

To get ready to differentiate this, we note $(x - y)/(1 - y) = 1 - (1 - x)/(1 - y)$ so we can write

$$x^* = g\left(\frac{1 - x}{1 - y}\right), \qquad \text{where } g(u) = \frac{(1 - u) + \sqrt{(1 - u)^2 + \sqrt{\mu} bu}}{2}.$$

To start to understand the properties of $g$ we note that

$$\begin{pmatrix} (x, y) & u = (1 - x)/(1 - y) & g(u) \\ (1, \sqrt{\mu}) & 0 & 1 \\ (\mu^{1/4}, \mu^{1/4}) & 1 & \mu^{1/4} \\ (\sqrt{\mu}, 1) & \infty & \sqrt{\mu} \end{pmatrix},$$

where for the last evaluation we note that $g(u) \approx \frac{1}{2}[(1 - u) + u - (1 - 2\sqrt{\mu})]$ for large $u$. Differentiating gives

$$\frac{\partial x^*}{\partial x} = g'\left(\frac{1 - x}{1 - y}\right)\frac{-1}{1 - y},$$

$$\frac{\partial x^*}{\partial y} = g'\left(\frac{1 - x}{1 - y}\right)\frac{1 - x}{(1 - y)^2},$$

$$\frac{\partial^2 x^*}{\partial x^2} = g''\left(\frac{1-x}{1-y}\right)\frac{1}{(1-y)^2},$$

$$\frac{\partial^2 x^*}{\partial y^2} = g''\left(\frac{1-x}{1-y}\right)\frac{(1-x)^2}{(1-y)^4} + g'\left(\frac{1-x}{1-y}\right)\frac{2(1-x)}{(1-y)^3},$$

where

$$g'(u) = \tfrac{1}{2}(-1 + \tfrac{1}{2}((1-u)^2 + 4\sqrt{\mu}u)^{-1/2}(-2(1-u) + 4\sqrt{\mu})),$$

$$g''(u) = -\tfrac{1}{8}((1-u)^2 + 4\sqrt{\mu}u)^{-3/2}(-2(1-u) + 4\sqrt{\mu})^2$$

$$+ \tfrac{1}{4}((1-u)^2 + 4\sqrt{\mu}u)^{-1/2}2.$$

**4.2.** *Itô's formula for the projection map.* Now consider the behavior of $(X_t^*, Y_t^*) = \Phi(X_t, Y_t)$. We conclude from Itô's formula that

$$X_t^* - X_0^* = \int_0^t g'\left(\frac{1-X_s}{1-Y_s}\right)\frac{-1}{1-Y_s}\,dX_s + \int_0^t g'\left(\frac{1-X_s}{1-Y_s}\right)\frac{1-X_s}{(1-Y_s)^2}\,dY_s$$

$$+ \frac{1}{2}\int_0^t g''\left(\frac{1-X_s}{1-Y_s}\right)\frac{1}{(1-Y_s)^2}X_s(1-X_s)\,ds$$

$$+ \frac{1}{2}\int_0^t \left[g''\left(\frac{1-X_s}{1-Y_s}\right)\frac{(1-X_s)^2}{(1-Y_s)^4} + g'\left(\frac{1-X_s}{1-Y_s}\right)\frac{2(1-X_s)}{(1-Y_s)^3}\right]$$

$$\times Y_s(1-Y_s)\,ds.$$

Using (2) to get expressions for $dX_s, dY_s$, the drift terms in the first two integrals which contain factors of $2N$ cancel leaving us with

$$\int_0^t g'\left(\frac{1-X_s}{1-Y_s}\right)\frac{-1}{1-Y_s}\sqrt{X_s(1-X_s)}\,dB_s^1$$

$$+ \int_0^t g'\left(\frac{1-X_s}{1-Y_s}\right)\frac{1-X_s}{(1-Y_s)^2}\sqrt{Y_s(1-Y_s)}\,dB_s^2$$

from which we see that the variance process is

$$\langle X^*\rangle_t - \langle X^*\rangle_0 = \int_0^t g'\left(\frac{1-X_s}{1-Y_s}\right)^2\frac{1}{(1-Y_s)^2}X_s(1-X_s)\,ds$$

$$+ \int_0^t g'\left(\frac{1-X_s}{1-Y_s}\right)^2\frac{(1-X_s)^2}{(1-Y_s)^4}Y_s(1-Y_s)\,ds.$$

For computing $E\tau$, it is convenient to have a symmetric diffusion so we will consider $X_t^* - Y_t^*$. Analogously to $x^*$ we have

$$y^* = \frac{(y-x) + \sqrt{(y-x)^2 + 4\sqrt{\mu}(1-x)(1-y)}}{2(1-x)} = g\left(\frac{1-y}{1-x}\right),$$

so

$$\frac{\partial y^*}{\partial y} = g'\left(\frac{1-y}{1-x}\right)\frac{-1}{1-x},$$

$$\frac{\partial y^*}{\partial x} = g'\left(\frac{1-y}{1-x}\right)\frac{1-y}{(1-x)^2},$$

$$\frac{\partial^2 y^*}{\partial y^2} = g''\left(\frac{1-y}{1-x}\right)\frac{1}{(1-x)^2},$$

$$\frac{\partial^2 y^*}{\partial x^2} = g''\left(\frac{1-y}{1-x}\right)\frac{(1-y)^2}{(1-x)^4} + g'\left(\frac{1-y}{1-x}\right)\frac{2(1-y)}{(1-x)^3}$$

and we conclude that $X_t^* - Y_t^*$ is a diffusion with two times the drift coefficient

$$
\begin{aligned}
(11) \quad & 2b(X_t, Y_t) \\
& = g''\left(\frac{1-X_t}{1-Y_t}\right)\frac{X_t(1-X_t)}{(1-Y_t)^2} - g''\left(\frac{1-Y_t}{1-X_t}\right)\frac{Y_t(1-Y_t)}{(1-X_t)^2} \\
& \quad + g''\left(\frac{1-X_t}{1-Y_t}\right)\frac{(1-X_t)^2 Y_t}{(1-Y_t)^3} + g'\left(\frac{1-X_t}{1-Y_t}\right)\frac{2(1-X_t)Y_t}{(1-Y_t)^2} \\
& \quad - g''\left(\frac{1-Y_t}{1-X_t}\right)\frac{(1-Y_t)^2 X_t}{(1-X_t)^3} - g'\left(\frac{1-Y_t}{1-X_t}\right)\frac{2(1-Y_t)X_t}{(1-X_t)^2}
\end{aligned}
$$

and diffusion coefficient

$$
\begin{aligned}
(12) \quad a(X_t, Y_t) &= g'\left(\frac{1-X_t}{1-Y_t}\right)^2\left[\frac{X_t(1-X_t)}{(1-Y_t)^2} + \frac{(1-X_t)^2 Y_t}{(1-Y_t)^3}\right] \\
& \quad + g'\left(\frac{1-Y_t}{1-X_t}\right)^2\left[\frac{Y_t(1-Y_t)}{(1-X_t)^2} + \frac{(1-Y_t)^2 X_t}{(1-X_t)^3}\right].
\end{aligned}
$$

4.3. *Staying near the curve.* To study the distance from the curve we will consider $\phi(X_t, Y_t) = (\mu - X_t^2 Y_t^2)^2$. It turns out that $\phi$ is a Lyapunov function for the system (3). In other words, $\phi(x, y) = (\mu - x^2 y^2)^2 \geq 0$, $\phi(x, y) = 0$ iff $(x, y)$ is a fixed point of the deterministic system (3), and for all $(x, y)$ in the neighborhood of fixed points $|\mu - x^2 y^2| \leq N^{-\delta}$ the change in the direction of the strong drift is

$$\nabla\phi \cdot F = -4xy(y(1-x) + x(1-y))\phi \leq -\beta\phi,$$

where $\beta = \inf_{x,y \, : \, |\mu - x^2 y^2| \leq N^{-\delta}}\{4xy(y(1-x) + x(1-y))\} > 0$.

Using Itô's formula we get

$$e^{2N\beta t}\phi(X_t, Y_t) = \int_0^t e^{2N\beta s}\left(2N\beta\phi(X_s, Y_s) + \frac{\partial\phi}{\partial x}2N(1-X_s)(\mu - X_s^2 Y_s^2)\right.$$

$$+ \frac{\partial \phi}{\partial y} 2N(1 - Y_s)(\mu - X_s^2 Y_s^2)$$

$$+ \frac{\partial^2 \phi}{\partial x^2} X_s(1 - X_s) + \frac{\partial^2 \phi}{\partial^2 y} \sqrt{Y_s(1 - Y_s)}\Big) ds$$

$$+ \int_0^t e^{2N\beta s} \Big( \frac{\partial \phi}{\partial x} \sqrt{X_s(1 - X_s)} \, dB_s^1 + \frac{\partial \phi}{\partial y} \sqrt{Y_s(1 - Y_s)} \, dB_s^2 \Big)$$

and collecting terms due to drift

$$\phi(X_t, Y_t) = e^{-2N\beta t} 2N \int_0^t e^{2N\beta s} \Big( \nabla \phi \cdot F + \beta \phi + \frac{\partial^2 \phi}{\partial x^2} X_s(1 - X_s)$$

$$+ \frac{\partial^2 \phi}{\partial^2 y} Y_s(1 - Y_s)\Big) ds$$

$$+ e^{-2N\beta t} \int_0^t e^{2N\beta s} \Big( \frac{\partial \phi}{\partial x} \sqrt{X_s(1 - X_s)} \, dB_s^1 + \frac{\partial \phi}{\partial y} \sqrt{Y_s(1 - Y_s)} \, dB_s^2 \Big).$$

For $(x, y)$ within a compact set it is easy to see that for any sequence of random times $\tau^N$ and $T \geq 0$ the second integral is less than

$$C \sup_{0 \leq t \leq T \wedge \tau^N} e^{-2N\beta t} \Big| \int_0^t e^{2N\beta s} (dB_s^1 + dB_s^2) \Big|,$$

which converges weakly to 0, likewise the third integral converges weakly to 0 as well.

Let $\tau^N = \inf\{t \geq 0 : |\mu - X_t^2 Y_t^2| > N^{-\delta}\}$, then for all $0 \leq t \leq \tau^N$ the first integral is negative, and since $\phi \geq 0$ this implies $\phi(X_{t \wedge \tau^N}, Y_{t \wedge \tau^N}) \Rightarrow 0$. In other words

$$(X_{t \wedge \tau^N}, Y_{t \wedge \tau^N}) - (X_{t \wedge \tau^N}^*, Y_{t \wedge \tau^N}^*) \quad \Rightarrow \quad 0.$$

4.4. *Limit theorem.* Our final step is to argue that as $N \to \infty$, $X_t^* - Y_t^*$ converges to the diffusion with coefficients $b(X_t^*, Y_t^*)$ and $a(X_t^*, Y_t^*)$. We have shown that for each $N$ the drift and diffusion coefficients are given by $b(X_t, Y_t)$ and $a(X_t, Y_t)$ as defined in (11) and (12). Since the coefficients $b$ and $a$ are bounded in a neighborhood of the curve of equilibria the sequence of processes is tight, and the weak convergence $(X_{t \wedge \tau^N}, Y_{t \wedge \tau^N}) - (X_{t \wedge \tau^N}^*, Y_{t \wedge \tau^N}^*) \Rightarrow 0$ that we have just shown together with Theorem 5.4 in Kurtz and Protter (1991) implies for all $0 \leq t \leq \tau^N$ we have convergence of $X_t^* - Y_t^*$ to the desired diffusion process.

**5. Proof of Theorem 3.** We follow the same general outline as the proof of Theorem 2, but there are some new steps which were trivial in the previous proof. First, finding a usable expression for the curve of equilibria requires more work, and to calculate the flow map requires the right change of variables to describe the surface along which it is constant. Next, we were not able to find a global Lyapunov function in this case, hence we linearize the system around the points on the curve, and investigate the behavior of the eigenvalues of the linearized system. As the last step we use the construction of Katzenberger to glue the local Lyapunov functions into a global one for the whole neighborhood of the curve.

5.1. *The curve of equilibria.* Consider the equations in (6) restricted to $x_1 = x_2 = x$, $y_1 = y_2 = y$.

$$(13) \quad \frac{dx_3}{dt} = -x_3 w + x_3 - 3bx_3 = x_3(\alpha - w),$$

$$(14) \quad \frac{dx}{dt} = -xw + x(y_3 + y) + bx_3 - 2bx = x(y_3 + y - 2b - w) + bx_3,$$

$$(15) \quad \frac{dy_3}{dt} = -y_3 w + y_3 - 3by_3 = y_3(\alpha - w),$$

$$(16) \quad \frac{dy}{dt} = -yw + y(x_3 + x) + by_3 - 2by = y(x_3 + x - 2b - w) + by_3,$$

where $w = x_3 + y_3 - x_3 y_3 + 2xy$ and $\alpha = 1 - 3b$. Letting $\beta = \alpha + 2b = 1 - b$, fixed points satisfy three equations

$$(17) \quad x_3 + y_3 - x_3 y_3 + 2xy = \alpha,$$

$$(18) \quad x(y_3 - \beta) + xy + bx_3 = 0,$$

$$(19) \quad y(x_3 - \beta) + xy + by_3 = 0.$$

To solve it is convenient to let

$$(20) \quad \Gamma = x_3 + y_3 - x_3 y_3 - \alpha = -2xy.$$

We can solve (17), (18), (19) for $x$ and $y$

$$(21) \quad x = \frac{\Gamma - 2bx_3}{2(y_3 - \beta)}, \qquad y = \frac{\Gamma - 2by_3}{2(x_3 - \beta)}.$$

Using these in (20) we have

$$(22) \quad 2\Gamma(x_3 - \beta)(y_3 - \beta) + (\Gamma - 2bx_3)(\Gamma - 2by_3) = 0.$$

Filling in the definition of $\Gamma$

$$2(-x_3 y_3 + x_3 + y_3 - \alpha)(x_3 y_3 - \beta x_3 - \beta y_3 + \beta^2)$$
$$+ (-x_3 y_3 + x_3(1 - 2b) + y_3 - \alpha)(-x_3 y_3 + x_3 + y_3(1 - 2b) - \alpha) = 0$$

and collecting terms gives

$$x_3^2 y_3^2 [-2+1] + (x_3^2 y_3 + x_3 y_3^2)[2\beta + 2b]$$
$$+ x_3 y_3 [-2\beta^2 - 4\beta + (1-2b)^2 + 1] + (x_3^2 + y_3^2)[-2\beta + (1-2b)]$$
$$+ (x_3 + y_3)[2\beta^2 + 2\alpha\beta - 2\alpha(1-b)] + [\alpha^2 - 2\alpha\beta^2] = 0,$$

which has the form $\sum_{i,j} c_{ij} x_3^i y_3^j$, where $c_{i,j} = c_{j,i}$. Using $\beta = 1-b$ and $\alpha = 1 - 3b$

$$c_{2,2} = -1,$$
$$c_{1,2} = c_{2,1} = 2\beta + 2b = 2,$$
$$c_{0,2} = c_{2,0} = -2\beta + (1-2b) = -1,$$
$$c_{1,1} = -2(1-2b+b^2) - 4(1-b) + (1-4b+4b^2) + 1 = -4 + 4b + 2b^2,$$
$$c_{0,1} = c_{1,0} = 2\beta^2 = 2 - 4b + 2b^2,$$
$$c_{0,0} = \alpha^2 - 2\alpha\beta^2 = \alpha(\alpha - 2\beta^2) = (1-3b)(-1+b-2b^2)$$
$$= -1 + 4b - 5b^2 + 6b^3.$$

To solve for $y_3$ note our equation has the form $d_0 + d_1 y_3 + d_2 y_3^2 = 0$ where $d_j = \sum_i c_{ij} x_3^i$, so

$$(23) \qquad y_3 = \frac{-d_1 \pm \sqrt{d_1^2 - 4d_0 d_2}}{2d_2}.$$

To see which root we want, note that if $x_3 = 0$ and $x = 0$ then $w = y_3 = \alpha$, so we want

$$y_3 = \frac{-c_{0,1} + \sqrt{c_{0,1}^2 - 4c_{0,0}c_{0,2}}}{c_{0,2}}$$

$$= \frac{-2\beta^2 + \sqrt{4\beta^4 - 4(\alpha^2 - 2\alpha\beta^2)(-1)}}{-2} = \alpha = 1 - 3b.$$

5.2. *Projection map.* It follows from the equations for the ODE that

$$\frac{d}{dt}\frac{y_3}{x_3} = \frac{-y_3}{x_3^2}\frac{dx_3}{dt} + \frac{1}{y_3}\frac{dy_3}{dt} = 0,$$

so $y_3/x_3$ is constant along solutions. This means that if $y_3/x_3 = r$ then the trajectory will flow the point on the curve of equilibria where $u = x_3^*$ has

$$(ru)^2 d_2(u) + (ru)d_1(u) + d_0(u) = 0.$$

Differentiating with respect to $r$, and letting $e_i = (u^i d_i)'$, where $'$ indicates the derivative with respect to $u$ we have

$$2ru^2 d_2 + r^2 e_2 \frac{du}{dr} + u d_1 + r e_1 \frac{du}{dr} + e_0 \frac{du}{dr} = 0.$$

Rearranging gives

$$(24) \qquad \frac{du}{dr} = -\frac{2ru^2 d_2 + u d_1}{r^2 e_2 + r e_1 + e_0}.$$

From this it follows that

$$\frac{d}{dx_3} u(y_3/x_3) = -\frac{y_3}{x_3^2} \frac{du}{dr} \frac{d}{dy_3} u(y_3/x_3) = \frac{1}{x_3} \frac{du}{dr},$$

$$\frac{d^2}{dx_3^2} u(y_3/x_3) = \frac{y_3^2}{x_3^4} \frac{d^2 u}{dr^2} + \frac{2y_3}{x_3^3} \frac{du}{dr} \frac{d^2}{dy_3^2} u(y_3/x_3) = \frac{1}{x_3^2} \frac{d^2 u}{dr^2}.$$

To now compute $d^2 u/dr^2$, we note that $du/dr = -f(r)/g(r)$ with $f'(r) = 2u^2 d_2 + (2re_2 + e_1) du/dr$ and $g'(r) = 2re_2 + e_1 + (r^2 e_2' + r e_1' + e_0')(du/dr)$ where $e_i'$ is the derivative of $e_i$ with respect to $r$.

Let $v = y_3^*$ and $q = x_3/y_3 = 1/r$. To compute the derivatives of $v$ we note that if $d_j(z) = \sum_i c_{ij} z^i$, then repeating the previous steps leads to

$$(25) \qquad \frac{dv}{dq} = -\frac{2qv^2 d_2(v) + v d_1(v)}{q^2 e_2(v) + q e_1(v) + e_0(v)}$$

and it follows that

$$\frac{d}{dx_3} v(x_3/y_3) = \frac{1}{y_3} \frac{dv}{dq}, \qquad \frac{d}{dy_3} v(x_3/y_3) = -\frac{x_3}{y_3^2} \frac{dv}{dq},$$

$$\frac{d^2}{dx_3^2} v(x_3/y_3) = \frac{1}{y_3^2} \frac{d^2 v}{dq^2}, \qquad \frac{d^2}{dy_3^2} v(x_3/y_3) = \frac{x_3^2}{y_3^4} \frac{d^2 v}{dq^2} + \frac{2x_3}{y_3^3} \frac{dv}{dq}.$$

To compute $d^2 v/dq^2$, we note that $dv/dq = -f_2(q)/g_2(q)$ with $f_2'(q) = 2v^2 d_2 + (2qe_2 + e_1) dv/dq$ and $g_2'(q) = 2qe_2 + e_1 + (q^2 e_2' + q e_1' + e_0')(dv/dq)$.

### 5.3. Itô's formula for the projection map.

Writing $R_s = Y_3(s)/X_3(s)$ we conclude from Itô's formula that

$$X_3^*(t) - X_3^*(0)$$

$$= \int_0^t u'(R_s) \cdot \frac{-Y_3(s)}{X_3^2(s)} \, dX_3(s) + \int_0^t u'(R_s) \cdot \frac{1}{X_3(s)} \, dY_3(s)$$

$$+ \frac{1}{2} \int_0^t \left[ u''(R_s) \cdot \frac{Y_3^2(s)}{X_3^4(s)} + u'(R_s) \cdot \frac{2Y_3(s)}{X_3^3(s)} \right] X_3(s)(1 - X_3(s)) \, ds$$

$$+ \frac{1}{2} \int_0^t u''(R_s) \cdot \frac{1}{X_3^2(s)} Y_3(s)(1 - Y_3(s)) \, ds.$$

The drift terms in the first two integrals cancel leaving us with

$$\int_0^t u'(R_s) \cdot \frac{-Y_3(s)}{X_3^2(s)} \sqrt{X_3(s)(1 - X_3(s))} \, dB_s^1$$

$$+ \int_0^t u'(R_s) \cdot \frac{1}{X_3(s)} \sqrt{Y_3(s)(1 - Y_3(s))} \, dB_s^2,$$

so the variance process

$$\langle X_3^* \rangle_t - \langle X_3^* \rangle_0$$
$$= \int_0^t u'(R_s)^2 \cdot \frac{Y_3^2(s)}{X_3^4(s)} X_3(s)(1 - X_3(s)) \, ds$$
$$+ \int_0^t u'(R_s)^2 \cdot \frac{1}{X_3(s)^2} Y_3(s)(1 - Y_3(s)) \, ds.$$

Let $Q_s = X_3(s)/Y_3(s)$. Combining these with the corresponding formulas for $Y_3^*(t)$ we see that $X_3^*(t) - Y_3^*(t)$ is a diffusion with two times the drift equal to

$$2b(X_3(t), Y_3(t))$$
$$= \left[ u''(R_t) \cdot \frac{Y_3^2(t)}{X_3^4(t)} + u'(R_t) \cdot \frac{2Y_3(t)}{X_3^3(t)} \right] X_3(t)(1 - X_3(t))$$
$$+ u''(R_t) \cdot \frac{1}{X_3^2(t)} Y_3(t)(1 - Y_3(t)) - v''(Q_t) \cdot \frac{1}{Y_3^2(t)} X_3(t)(1 - X_3(t))$$
$$- \left[ v''(Q_t) \cdot \frac{X_3^2(t)}{Y_3^4(t)} + v'(Q_t) \cdot \frac{2X_3(t)}{Y_3^3(t)} \right] Y_3(t)(1 - Y_3(t))$$

and variance

$$a(X_3(t), Y_3(t))$$
$$= u'(R_t)^2 \cdot \left[ \frac{Y_3^2(t)}{X_3^4(t)} X_3(t)(1 - X_3(t)) + \frac{1}{X_3(t)^2} Y_3(t)(1 - Y_3(t)) \right]$$
$$+ v'(Q_t)^2 \cdot \left[ \frac{1}{Y_3(t)^2} X_3(t)(1 - X_3(t)) + \frac{X_3^2(t)}{Y_3^4(t)} Y_3(t)(1 - Y_3(t)) \right].$$

5.4. *Linearization.* In order to show that the diffusion stays close to the curve of equilibria, we will investigate the linearization of the ODE around these fixed points. Recall from (6) the equations are

$$dx_3/dt = -x_3 w + x_3 - 3bx_3,$$
$$dx_2/dt = -x_2 w + x_2(y_3 + y_1) + bx_3 - 2bx_2,$$
$$dx_1/dt = -x_1 w + x_1(y_3 + y_2) + bx_3 - 2bx_1,$$
$$dy_3/dt = -y_3 w + y_3 - 3by_3,$$
$$dy_2/dt = -y_2 w + y_2(x_3 + x_1) + by_3 - 2by_2,$$
$$dy_1/dt = -y_1 w + y_1(x_3 + x_2) + by_3 - 2by_1$$

with $w = x_3 + y_3 - x_3 y_3 + x_2 y_1 + x_1 y_2$.

We write the coordinates $z = (z_1, z_2, z_3, z_4, z_5, z_6) = (x_3, y_3, x_2, x_1, y_2, y_1)$. If we have a dynamical system of the form $dz_i/dt = F_i(z)$ and we let $z = z^* + \varepsilon p$ where $z^*$ is a fixed point, then comparing terms of order $\varepsilon$ we see that the linearization at the fixed point is

$$\frac{dp_i}{dt} = \sum_j \frac{\partial F_i}{\partial z_j}(z^*) p_j.$$

To compute the linearization it is useful to note that $\nabla w = (1 - y_3, 1 - x_3, y_1, y_2, x_1, x_2)$ and that at $z^*$ we have $w = \alpha = 1 - 3b$

$$\nabla F_1(z^*) = -x_3 \nabla w,$$

$$\nabla F_2(z^*) = -y_3 \nabla w,$$

$$\nabla F_3(z^*) = -x_2 \nabla w + (b, x_2, y_3 + y_1 - 1 + b, 0, 0, x_2),$$

$$\nabla F_4(z^*) = -x_1 \nabla w + (b, x_1, 0, y_3 + y_2 - 1 + b, x_1, 0),$$

$$\nabla F_5(z^*) = -y_2 \nabla w + (y_2, b, 0, y_2, x_3 + x_1 - 1 + b, 0),$$

$$\nabla F_6(z^*) = -y_1 \nabla w + (y_1, b, y_2, 0, 0, x_3 + x_2 - 1 + b).$$

On the curve of equilibria $x_2 = x_1 = x$ and $y_2 = y_1 = y$, so the four-dimensional space $p_3 = p_4$, $p_5 = p_6$ is invariant for the linearization, as is the orthogonal two-dimensional subspace of vectors of the form: $(0, 0, u, -u, v, -v)$. Thus our six-dimensional problem decomposes into a two-dimensional one and a four-dimensional one. Taking the smaller one first, multiplying the linearization on the right the linear map in this subspace is

$$(u, v) \to (u(y_3 + y - 2b - w) - vx, -uy + v(x_3 + x - 2b - w)).$$

To simplify we note that in equilibrium $y_3 + y - 2b - w = -bx_3/x$ and $x_3 + x - 2b - w = -by_3/y$ so the matrix is

$$\begin{pmatrix} -bx_3/x & -x \\ -y & -by_3/y \end{pmatrix}.$$

For a $2 \times 2$ matrix $M$, $\text{trace}(M) = \lambda_1 + \lambda_2$ and $\det(M) = \lambda_i \lambda_2$, so the real part $\Re(\lambda_i) < 0$ if and only if $\text{trace}(M) < 0$ and $\det(M) > 0$. Clearly $\text{trace}(M) = -bx_3/x - by_3/y < 0$. The sign of $\det(M) = b^2 x_3 y_3/xy - xy$ is not so clear but as it turns out (for the proof see Appendix).

LEMMA 2. $bx_3/x > y$ and $by_3/y > x$ so $b^2 x_3 y_3/xy > xy$.

Turning to the four-dimensional subspace, we use coordinates $(z_1, z_2, z_3, z_4) = (x_3, y_3, x, y)$, and the linearization reduces to

$$\mathcal{F} = \begin{pmatrix} -x_3(1 - y_3) & -x_3(1 - x_3) & -2x_3 y & -2x_3 x \\ -y_3(1 - y_3) & -y_3(1 - x_3) & -2y_3 y & -2y_3 x \\ -x(1 - y_3) + b & -x(1 - x_3) + x & -2xy - bx_3/x & -2x^2 + x \\ -y(1 - y_3) + y & -y(1 - x_3) + b & -2y^2 + y & -2xy - by_3/y \end{pmatrix},$$

where we have again used $y_3 + y - 2b - w = -bx_3/x$ and $x_3 + x - 2b - w = -by_3/y$.

To simplify, we will change to a new basis $v_1$, $v_2$, $v_3$, $v_4$. To see how this will affect the linearization, note that if $V$ is the matrix with $i$th column $v_i$ then a vector with coordinates $q$ in terms of the new basis, is in the original coordinate system

$$p_j = \sum_k v_{j,k} q_k = Vq.$$

To return to the $q$ system $q = V^{-1}p$, so

$$\frac{dq}{dt} = V^{-1}\mathcal{F}Vq.$$

To cancel the multiples of $\nabla w$ that appear in all of the rows while preserving the symmetry in the problem, we choose

$$V^{-1} = \begin{pmatrix} 1/x_3 & -1/y_3 & 0 & 0 \\ 1/x_3 & 1/y_3 & 0 & 0 \\ -x/x_3 & 0 & 1 & 0 \\ 0 & -y/y_3 & 0 & 1 \end{pmatrix}, \qquad V = \begin{pmatrix} x_3/2 & x_3/2 & 0 & 0 \\ -y_3/2 & y_3/2 & 0 & 0 \\ x/2 & x/2 & 1 & 0 \\ -y/2 & y/2 & 0 & 1 \end{pmatrix}.$$

It is easy to check that $V^{-1}V = I$. The left multiplication simplifies

$$V^{-1}\mathcal{F} = \begin{pmatrix} 0 & 0 & 0 & 0 \\ -2(1-y_3) & -2(1-x_3) & -4y & -4x \\ b & x & -bx_3/x & x \\ y & b & y & -by_3/y \end{pmatrix},$$

while the right multiplication make things a little messier,

$$V^{-1}\mathcal{F}V$$
$$= \begin{pmatrix} 0 & 0 & 0 & 0 \\ -x_3(1-y_3)+y_3(1-x_3) & -x_3(1-y_3)-y_3(1-x_3)-4xy & -4y & -4x \\ -x(y_3+y)/2 & x(y_3+y)/2 & -bx_3/x & x \\ y(x_3+x)/2 & y(x_3+x)/2 & y & -by_3/y \end{pmatrix}.$$

This matrix has one zero eigenvector corresponding to the direction along the curve. To find conditions which guarantee that all of the eigenvalues of the three-dimensional operator on the perpendicular subspace have negative real parts, we turn to the Appendix 2 of Murray (1989). Given the characteristic polynomial

$$P(\lambda) = \lambda^n + a_1\lambda^{n-1} + \cdots + a_n = 0$$

the *Routh Hurwitz conditions* in terms of the $a$'s are necessary and sufficient. For a $3 \times 3$ matrix $M = [m_{ij}]$ these are

$$a_1 > 0, \qquad a_3 > 0, \qquad a_1a_2 - a_3 > 0,$$

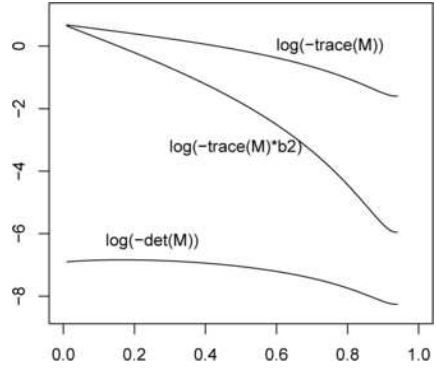

Fig. 7. *Picture of the of coefficients determining the behavior of the linearization of the deterministic drift at the curve of equilibria in the subfunctionalization model.*

which in terms of $b_1 = \text{trace}(M)$, $b_3 = \det(M)$, and

$$b_2 = m_{22}m_{33} + m_{11}m_{33} + m_{11}m_{22} - m_{23}m_{32} - m_{12}m_{21} - m_{13}m_{31}$$

$$= (m_{11}m_{22} - m_{12}m_{21}) + (m_{22}m_{33} - m_{23}m_{32}) + (m_{11}m_{33} - m_{13}m_{31})$$

are

$$(26) \qquad \text{trace}(M) < 0, \qquad \det(M) < 0, \qquad \det(M) > \text{trace}(M)b_2.$$

Note that for this to be possible we must have $b_2 > 0$.

Before trying to prove these three conditions, we computed them numerically for the concrete example $b = 0.001$. In Figure 7 we have plotted the values. Since the determinant turns out to be small, we have plotted the logarithm of $-\text{trace}(M)$, $-\det(M)$ and $-\text{trace}(M)b_2$.

The eigenvalues of our $4 \times 4$ matrix $V^{-1}\mathcal{F}V$ are a 0 and the eigenvalues of the $3 \times 3$ matrix

$$(27) \qquad M = \begin{pmatrix} -x_3(1-y_3) - y_3(1-x_3) - 4xy & -4y & -4x \\ x(y_3+y)/2 & -bx_3/x & x \\ y(x_3+x)/2 & y & -by_3/y \end{pmatrix}.$$

It is clear that $\text{trace}(M) < 0$. The determinant is

$$\det(M) = (-x_3(1-y_3) - y_3(1-x_3) - 4xy)\left[\frac{b^2 x_3 y_3}{xy} - xy\right]$$

$$- 4y\left[\frac{x(y_3+y)}{2}\left(\frac{by_3}{y}\right) + \frac{xy(x_3+x)}{2}\right]$$

$$- 4x\left[\frac{xy(y_3+y)}{2} + \left(\frac{bx_3}{x}\right)\frac{y(x_3+x)}{2}\right],$$

so Lemma 2 implies $\det(M) < 0$.

For the third condition, since $m_{11} < 0$ and since Lemma 2 implies that $b^2 x_3 y_3 / xy > xy$, we have

$$
\begin{aligned}
(28) \quad b_2 = & -\frac{bx_3}{x} m_{11} + 2xy(y_3 + y) + \left(\frac{b^2 x_3 y_3}{xy} - xy\right) \\
& -\frac{by_3}{y} m_{11} + 2xy(x_3 + x) > 0,
\end{aligned}
$$

which is necessary for $\det(M) > b_2 \operatorname{trace}(M)$. Once

LEMMA 3.   $\det(M) > b_2 \operatorname{trace}(M)$.

Lemma 3 is proved (see the Appendix), we will have verified the third condition in (26) and shown that all the eigenvalues have negative real part.

5.5. *Staying near the curve.*   In order to construct a local Lyapunov function in a neighborhood of the whole curve of equilibria with exponential convergence properties we follow closely the arguments of Katzenberger. We present the outline of the proof and refer the reader to Section 3 of Katzenberger (1991) for more details.

To begin we consider a linear dynamical system

$$
(29) \qquad \frac{dx(t)}{dt} = Ax(t).
$$

If things are nice, for example, if all eigenvalues are distinct, when we let $V$ be the matrix with columns equal to the eigenvectors of $A$, we will have

$$
A = V\Lambda V^{-1},
$$

where $V$ is a diagonal matrix with $\Lambda_{ii} = \lambda_i$ the eigenvalues of $A$. Let $e^{tA} = \sum_{k=0}^{\infty} \frac{(tA)^k}{k!}$. The solution of (29) can be written

$$
x(t) = e^{tA} x(0) = V e^{\Lambda t} V^{-1} x(0).
$$

Using a superscript $T$ for transpose, a little linear algebra shows

$$
\|V^{-1} x(t)\|^2 = x(0)^T (V^{-1})^T e^{2\Lambda t} V^{-1} x(0) = \sum_i e^{2\lambda_i t} (V^{-1} x(0))_i^2,
$$

so if all eigenvalues have negative real part $\|V^{-1} x(t)\|$ decreases exponentially.

To do this in general for a nonlinear dynamical system

$$
\frac{dx(t)}{dt} = F(x(t))
$$

consider the linearized system with $A = [\frac{\partial F_i}{\partial z_j}(z^*)]_{i,j}$ along the curve of fixed points $z^*$. If all of the eigenvalues of $A$ have negative real part, then there is a $\delta > 0$ so that they all lie in $D(\delta) = \{z \in \mathbf{C} : |\delta z + 1| < 1\}$. $B = I + \delta A$ has spectral radius $\rho(B) < 1$. $\rho(B) = \inf_{Q \in Gl(d)} |Q^{-1}BQ|$ so there is an invertible matrix $V$ with

$$\|V^{-1}BV\| < 1.$$

Since $A = (B - I)/\delta$ and $I$ and $B$ commute, we have

$$e^{tA} = e^{Bt/\delta} e^{-t/\delta}$$

and a little linear algebra shows that $\|V^{-1}x(t)\|$ decreases exponentially at rate $\alpha = 1 - \rho(B)$.

We now construct a Lyapunov function in the neighborhood of an arbitrary fixed point $x$ on the curve. Consider first the case that the curve of fixed points in the neighborhood $W$ of any fixed point $x$ is of form $N \cap W$ for some linear subspace $N$. Then $N = \mathrm{Ker}(A)$ and $P = \mathrm{Ran}(A)$ decompose the whole space, and the spectral radius of $B = I + A$ restricted to $P$ is $\rho(B)|_P < 1$. There is an invertible matrix $V$ on $P$ with $\|V^{-1}B|_P V\| < 1$. Extend $V$ to the whole space by setting $V = V\Pi_P + \Pi_N$ with projection maps $\Pi_P$, and $\Pi_N$ onto $P$ and $N$, respectively. Then $\|V^{-1}\Pi_P x(t)\|$ has an exponential decrease with rate $\alpha = 1 - \rho(B)|_P$.

Let $\theta(y) = F(y) - Ay$, then $\mathrm{D}\theta(x) = [\frac{\partial \theta_i}{\partial x_j}]_{i,j} \equiv 0$, and for $\varepsilon > 0$ there is a neighborhood $W' \subset W$ of $x$ such that $\|V^{-1}\mathrm{D}\theta(y)V\| \le \varepsilon$ for all $y \in W'$, and consequently $\|V^{-1}\theta(x)\| \le \varepsilon \|V^{-1}\Pi_P x\|$. Take $\varepsilon < \alpha$, and let $\tau = \inf\{t \ge 0 : x(t) \notin W'\}$. Integration by parts on $e^{-At}x(t)$ gives

$$x(t) = e^{A(t-s)}x(s) + \int_s^t e^{A(t-r)}\theta(x(r))\, dr$$

hence

$$\|V^{-1}\Pi_P x(t)\| \le e^{-\alpha(t-s)}\|V^{-1}\Pi_P x(s)\| + \int_s^t e^{-\alpha(t-r)}\|V^{-1}\Pi_P \theta(x(r))\|\, dr$$

$$\le e^{-\alpha(t-s)}\|V^{-1}\Pi_P x(s)\| + \varepsilon \int_s^t e^{-\alpha(t-r)}\|V^{-1}\Pi_P x(r)\|\, dr.$$

Gronwall's inequality implies

$$(30) \qquad \|V^{-1}\Pi_P x(t)\| \le \|V^{-1}\Pi_P x(0)\| e^{-(\alpha-\varepsilon)t}.$$

Since $\Pi_N A = 0$ and $\Pi_N e^A = \Pi_N$,

$$(31) \qquad \begin{aligned}
&\|V^{-1}\Pi_N x(t) - V^{-1}\Pi_N x(0)\| \\
&= \int_0^t \|V^{-1}\Pi_N \theta(x(r))\, dr\| \\
&\le \varepsilon \int_0^t \|V^{-1}\Pi_N x(r)\|\, dr \le \varepsilon \|V^{-1}\Pi_P x(0)\| \int_0^t e^{-(\alpha-\varepsilon)r}\, dr \\
&= \frac{\varepsilon}{\alpha-\varepsilon}\|V^{-1}\Pi_P x(0)\|(1 - e^{-(\alpha-\varepsilon)t}).
\end{aligned}$$

Adding $\varepsilon/(\alpha - \varepsilon)$ times (30) to (31)

$$\|V^{-1}\Pi_N x(t)\| + \frac{\varepsilon}{\alpha - \varepsilon}\|V^{-1}\Pi_P x(t)\| \le \|V^{-1}\Pi_N x(0)\| + \frac{\varepsilon}{\alpha - \varepsilon}\|V^{-1}\Pi_P x(0)\|.$$

So, if $\rho > 0$ is such that $W'' = \{x : \|V^{-1}\Pi_N x\| + \frac{\varepsilon}{\alpha-\varepsilon}\|V^{-1}\Pi_P x\| < \rho\} \subset W'$, then for any starting point in $W''$ the exit time $\tau = \infty$. Let

$$f(x) = \frac{\|V^{-1}\Pi_N x\|}{\rho}, \qquad\qquad g(x) = \frac{\varepsilon}{\alpha - \varepsilon}\frac{\|V^{-1}\Pi_P x\|}{\rho},$$

$$v(x) = g(x) + \sqrt{f^2(x) + g^2(x)}, \qquad u(x) = (1 - v(x))^3$$

and $W''' = \{x \in W'' : u(x) > 0\}$. On $W'''$ let $\phi(x) = g^2(x)/u(x)$. Then,

$$v(x(t)) \le v(x(0)), \qquad \phi(x(t)) \le e^{-2(\alpha-\varepsilon)t}\phi(x(0)), \qquad \nabla\phi \cdot F \le -2(\alpha - \varepsilon)\phi.$$

Note that $\phi$ satisfies $\phi \ge 0$ $\phi(x) = 0$ if $x$ is a fixed point, and the reason for using the function $u(x)$ in the construction is that $\phi(x) = \infty$ for $x$ outside $W'''$.

In general let $N$ be the tangent space of the curve of fixed points at $x$. Then $N = \mathrm{Ker}(A)$ and $P = \mathrm{Ran}(A)$ decompose the full space and there is a smooth coordinate function $\eta$ that will map any $x \in N$ to $\eta(x) \in P$ in such a way that $x + \eta(x)$ is a fixed point. Using this coordinate transformation we can construct a Lyapunov function $\phi$ in a neighborhood of $x$ in the analogous fashion to the one above.

We can now patch these local Lyapunov functions together over the whole curve of equilibria. If $\phi_1$ is a Lyapunov function constructed for the neighborhood $W_1$ and $\phi_2$ for $W_2$, define a function on $W = W_1 \cup W_2$ by $\phi(x) = 0$ if $x$ is a fixed point and $\phi(x) = 1/(\phi_1(x) + \phi_2(x))$ for any other point of $W$. It is easy to verify $\phi$ has the same properties as the functions $\phi_i$.

Given a local exponential Lyapunov function $\phi$ the proof that the diffusion $Z(t) = (X_3(t), Y_3(t), X_2(t), X_1(t), Y_2(t), Y_1(t))$ converges weakly to the curve of equilibria and the proof that the limiting process is characterized by an SDE derived for $\Phi(Z_t)$ follow by the same arguments as those described in Sections 3.3. and 3.4.

## APPENDIX

PROOF OF LEMMA 2. In equilibrium $w = \alpha = 1 - 3b$, so by (13)

$$0 < bx_3/x = 1 - b - y_3 - y, \qquad 0 < by_3/y = 1 - b - x_3 - x.$$

Thus if can show $y_3 + y_2 + y_1 < 1 - b$ (and hence by symmetry $x_3 + x_2 + x_1 < 1 - b$) the desired result will follow. To do this we note that $w = \alpha = 1 - 3b$ so

$$\frac{dy_0}{dt} = by_3 + 2by_2 + 2by_1 + x_3 y_0 - y_0(1 - 3b).$$

Thus in equilibrium $b(1 - y_0) + b(y_1 + y_2) = y_0(1 - 3b - x_3)$ and it follows that $y_0 \geq b/(1 - 2b - x_3) > b$. $\quad \square$

PROOF OF LEMMA 3. To check the condition $\det(M) > b_2 \operatorname{trace}(M)$ we rearrange the formula (28) for the determinant of $M = [m_{ij}]$ (27)

$$\det(M) = m_{11}\left[\frac{b^2 x_3 y_3}{xy} - xy\right] - \frac{by_3}{y} 2xy(y_3 + y) - \frac{bx_3}{x} 2xy(x_3 + x)$$
$$- y \cdot 2xy(x_3 + x) - x \cdot 2xy(y_3 + y).$$

The trace of $M$, consists of three negative terms $m_{ii}$, while $b_2$ consists of five positive terms which we will refer to as $b_{2j}$ where $j$ is the order in (28). To prove the desired inequality it is enough to find terms smaller (i.e., more negative) than the five parts of $\det(M)$, $d_1, \ldots, d_5$, within the fifteen products in $b_2 \operatorname{trace}(M)$. The first three are easy: $d_1 = m_{11}b_{23}$, $d_2 = m_{33}b_{22}$, $d_3 = m_{22}b_{25}$. For the final two we will use the next lemma.

LEMMA A.1. *If $b \leq 0.01$ then in an equilibrium we have $y_3(1 - x_3) \geq x \geq y$ whenever $x_3 \geq y_3$.*

Using the result of Lemma A.1 we see that when $x_3 \geq y_3$, $m_{11}(b_{22} + b_{25}) < -y_3(1 - x_3)(b_{22} + b_{25}) \leq d_4 + d_5$. Symmetry implies that when $x_3 \leq y_3$, $m_{11}(b_{22} + b_{25}) < -x_3(1 - y_3)(b_{22} + b_{25}) \leq d_4 + d_5$ as well. $\quad \square$

PROOF OF LEMMA A.1. We want to show that if $b \leq 0.01$ then in equilibrium,

$$(32) \qquad\qquad y_3(1 - x_3) \geq x \geq y,$$

whenever $x_3 \leq y_3$. Recall that if $\beta = 1 - b$ and $\alpha = 1 - 3b$, the fixed points satisfy three equations

$$(33) \qquad\qquad x_3 + y_3 - x_3 y_3 + 2xy = \alpha,$$

$$(34) \qquad\qquad x(y_3 - \beta) + xy + bx_3 = 0,$$

$$(35) \qquad\qquad y(x_3 - \beta) + xy + by_3 = 0.$$

For the second inequality we note that taking the difference (34) minus (35):

$$x(y_3 - \beta) - y(x_3 - \beta) = b(y_3 - x_3).$$

Since $x_3 \leq y_3 < \beta = 1 - b$ (recall that in the proof of Lemma 2 we showed $y_0 > b$) it follows that

$$y = \left(\frac{\beta - y_3}{\beta - x_3}\right)x - \frac{b(y_3 - x_3)}{\beta - x_3} \leq x.$$

To prove the first inequality in (32) we will show that $y_3(1-x_3) - x$ is a decreasing function $h(x_3)$ plus an error of $O(b)$, and the value of $h$ at the point where $x_3 = y_3$ is positive and $O(\sqrt{b})$. Our equation (23) for $y_3$ in terms of $x_3$ has the form

$$y_3 = \frac{-d_1 + \sqrt{d_1^2 - 4d_0 d_2}}{2d_2},$$

where $d_j = \sum_i c_{ij} x_3^i$, and in matrix form

$$c_{i,j} = \begin{pmatrix} -1 + 4b - 5b^2 + b^3 & 2 - 4b + 2b^2 & -1 \\ 2 - 4b + 2b^2 & -4 + 4b + 2b^2 & 2 \\ -1 & 2 & -1 \end{pmatrix}$$

with the rows and columns in the order $0, 1, 2$.

$$d_2 = (-1 + 2x_3 - x_3^2) = -(1-x_3)^2,$$

$$d_1 = (2 - 4x_3 + 2x_3^2) + b(-4 + 4x_3) + 2b^2(1+x_3)$$

$$= 2(1-x_3)^2 - 4b(1-x_3) + 2b^2(1+x_3),$$

$$d_0 = (-1 + 2x_3 - x_3^2) + b(4 - 4x_3) - b^2(5 - 2x_3) + b^3$$

$$= -(1-x_3)^2 + 4b(1-x_3) - b^2(5 - 2x_3) + b^3.$$

If $b = 0$ then $d_0 + d_2 = -d_1$ so

$$(1 - y_3)(d_0 - d_2 y_3) = d_0 + d_1 y_3 + d_2 y_3^2$$

and $y_3 = 1$. One can also see this from

$$d_1^2 - 4d_0 d_2 = d_1^2 + 4(d_1 + d_2)d_2 = (d_1 + 2d_2)^2 = 0.$$

Using the definitions of the $d_i$, we have

$$d_1^2 = 4(1-x_3)^4 - 16b(1-x_3)^3 + 16b^2(1-x_3)^2 + 8b^2(1-x_3)^2(1+x_3)$$

$$- 16b^3(1-x_3)(1+x_3) + 4b^4(1+x_3)^2,$$

$$-4d_2 d_0 = -4(1-x_3)^4 + 16b(1-x_3)^3 - 4b^2(1-x_3)^2(5-2x_3) + 4b^3(1-x_3)^2.$$

Adding the last two equations gives

$$d_1^2 - 4d_0 d_2 = 16b^2(1-x_3)^2 + 8b^2(1-x_3)^2(1+x_3) - 4b^2(1-x_3)^2(5-2x_3)$$

$$- 16b^3(1-x_3)(1+x_3) + 4b^3(1-x_3)(1-x_3) + 4b^4(1+x_3)^2$$

$$= 4b^2(1-x_3)^2[1 + 4x_3] - 4b^3(1-x_3)[3 + 5x_3] + 4b^4(1+x_3)^2$$

$$= 4b^2(1-x_3)^2[1 + 4x_3] \cdot M,$$

where the mess

$$M = 1 - \frac{b}{(1-x_3)} \cdot \frac{3 + 5x_3}{(1 + 4x_3)} + \frac{b^2}{(1-x_3)^2} \cdot \frac{(1+x_3)^2}{(1 + 4x_3)}.$$

Using $2d_2y_3 = -d_1 + \sqrt{d_1^2 - 4d_0d_2}$, we see that

$$-2(1-x_3)^2 y_3 = -2(1-x_3)^2 + 4b(1-x_3) - 2b^2(1+x_3)$$
$$+ 2b(1-x_3)\sqrt{1+4x_3}M^{1/2}$$

and we have

$$y_3 = 1 - \frac{b}{1-x_3}(2 + \sqrt{1+4x_3})$$
$$+ \frac{b^2}{(1-x_3)^2}\left(1 + x_3 - (1-x_3)\sqrt{1+4x_3}\frac{M^{1/2}-1}{b}\right)$$
$$\equiv 1 - q(x_3)b + r(x_3, b)b^2.$$

Note that $q(x_3) \to 3$ as $x_3 \to 0$, which is consistent with our earlier calculation that $y_3 \to 1 - 3b$ and we have $q(x_3) \geq 3$ for $x_3 \in [0,1]$.

To find $x$ using (21), we first compute

$$\Gamma = y_3 + x_3(1-y_3) - (1-3b)$$
$$= -qb + rb^2 + x_3(qb - rb^2) + 3b$$
$$= b(3 - q(1-x_3)) + rb^2(1-x_3)$$

and then

$$x = \frac{\Gamma - 2bx_3}{2(y_3 - \beta)} = \frac{b(3 - q(1-x_3)) + rb^2(1-x_3) - 2bx_3}{2[1 - qb + rb^2 - (1-b)]}.$$

Multiplying top and bottom by $-1/b$

$$(36) \qquad x = \frac{q(1-x_3) + (2x_3 - 3) + br(1-x_3)}{2(q-1) - 2rb}.$$

To check the condition $y_3(1-x_3) > x$ now we note that ignoring terms of order $b$:

$$y_3(1-x_3) - x = 1 - x_3 - \frac{q(1-x_3) + (2x_3 - 3)}{2(q-1)}$$
$$= 1 - x_3 - \frac{2 + \sqrt{1+4x_3} + 2x_3 - 3}{2(1 + x_3 + \sqrt{1+4x_3})} \cdot (1-x_3)$$
$$= (1-x_3)\frac{3 + \sqrt{1+4x_3}}{2(1 + x_3 + \sqrt{1+4x_3})} \equiv h(x_3).$$

It is not immediately obvious that $h$ is decreasing but writing $h = (1-x_3)f/(2g)$, where $f, g \geq 0$ on $[0,1]$ we have

$$\frac{dh}{dx_3} = -\frac{f}{2g} + \frac{(1-x_3)}{2} \cdot \frac{(4/2)(1+4x_3)^{-1/2}g - (1 + (4/2)(1+4x_3)^{-1/2})f}{g^2}.$$

The numerator of the final fraction is

$$2(1 + 4x_3)^{-1/2}[1 + x_3 + \sqrt{1 + 4x_3}] - [1 + 2(1 + 4x_3)^{-1/2}](3 + \sqrt{1 + 4x_3})$$

$$= (2 + 2x_3 - 6)(1 + 4x_3)^{-1/2} + 2 - 2 - 3 - \sqrt{1 + 4x_3} < 0,$$

so $dh/dx_3 < 0$. To complete the proof now, we will first approximately evaluate $h(x_3)$ at the point where $x_3 = y_3$, and then investigate the errors in our approximations.

*Solution when* $x_3 = y_3$. We proceed by guessing and then verifying our answer. From the previous calculation $y_3 = 1 - Qb/(1 - x_3)$, where $Q \approx 2 + \sqrt{5}$. If $x_3 = 1 - u\sqrt{b} = y_3$, then

$$u = \frac{2 + \sqrt{5}}{u} \quad \text{and hence} \quad u = \sqrt{2 + \sqrt{5}}.$$

From (36), since $q$ is large, $x_3 \approx 1$, and $q = u/\sqrt{b}$

$$x \approx \frac{q(1 - x_3) - 1}{2q} \approx \frac{2 + \sqrt{5} - 1}{2u/\sqrt{b}},$$

so we have $x = y = v\sqrt{b}$ with $v = \frac{1 + \sqrt{5}}{2\sqrt{2 + \sqrt{5}}}$.

To derive this from the equations for the equilibria, we note that when $x_3 = y_3$ and $x = y$, (33) implies:

$$1 - (1 - x_3)^2 + 2x^2 = 1 - 3b,$$

$$-u^2b + 2v^2b = -3b$$

or $u^2 - 2v^2 = 3$. Using (34) and dropping terms of order $b^{3/2}$

$$x(x_3 - (1 - b)) + x^2 + bx_3 = 0,$$

$$-vub + v^2b + b = 0.$$

The solution to the quadratic $v^2 - uv + 1 = 0$ is

$$v = \frac{u \pm \sqrt{u^2 - 4}}{2}.$$

Taking $u = \sqrt{2 + \sqrt{5}}$ and choosing the minus root gives

$$v = \frac{\sqrt{2 + \sqrt{5}} - \sqrt{\sqrt{5} - 2}}{2}$$

$$= \frac{2 + \sqrt{5} - \sqrt{(\sqrt{5} - 2)(2 + \sqrt{5})}}{2\sqrt{2 + \sqrt{5}}} = \frac{1 + \sqrt{5}}{2\sqrt{2 + \sqrt{5}}}.$$

For the other equation we note that

$$u^2 - 2v^2 = 2 + \sqrt{5} - 2\frac{6 + 2\sqrt{5}}{4(2 + \sqrt{5})} = \frac{9 + 4\sqrt{5} - 3 - \sqrt{5}}{2 + \sqrt{5}} = 3.$$

At $x_3 = y_3$ we have

$$y_3(1 - x_3) - x \approx u\sqrt{b} - \frac{u^2 - 1}{2u}\sqrt{b} = \frac{u^2 + 1}{2u}\sqrt{b} > 0.$$

*Error analysis.* To evaluate the error, we write $q = Q/(1 - x_3)$ and $r = R/(1 - x_3)^2$. Taking the difference between the exact expression in (36) and the approximation obtained by setting $b = 0$ gives

$$= \frac{q(1 - x_3) + (2x_3 - 3) + br(1 - x_3)}{2(q - 1) - 2rb} = \frac{q(1 - x_3) + (2x_3 - 3)}{2(q - 1)}$$

$$= \frac{br[-(1 - x_3)2(q - 1) + q(1 - x_3) + (2x_3 - 3)]}{4(q - 1 - rb)(q - 1)}$$

$$= \frac{bR[-Q - 1]}{4(Q - 1 - Rb/(1 - x_3))(Q - 1)}.$$

We have $Q + 1 \leq 3 + \sqrt{5}$. To bound $R$ we begin by noting that if $b \leq 0.01$ and $1 - x_3 \geq b^{1/2}$ then

$$1 \geq M \geq 1 - 3b/(1 - x_3) \quad \text{and hence} \quad 0 \leq (1 - x_3)\frac{1 - M^{1/2}}{b} \leq 3,$$

which implies $0 \leq R \leq 2 + 3\sqrt{5} = 8.7082$. Since $R\sqrt{b} \leq 1$ it follows that when $1 - x_3 \geq b^{1/2}$ the absolute value of the error is

$$\leq b\frac{(2 + 3\sqrt{5})(3 + \sqrt{5})}{8} = \frac{19 + 7\sqrt{5}}{8} = 4.331b.$$

To see if this is small enough, we note that $u = \sqrt{2 + \sqrt{5}} = 2.05817$ so $(u^2 - 1)/2u = 1.27202$, and we do have $1.272\sqrt{b} > 3.906b$ when $b < 0.01$. The last detail to check is that when $x_3 = b^{1/2}$ and $b \leq 0.01$

$$y_3 \leq 1 - 3b^{1/2} + 8.7082b < 1 - 2b^{1/2},$$

so the assumption $1 - x_3 \geq b^{1/2}$ we used to bound $M$ is valid. □

DEPARTMENT OF MATHEMATICS
CORNELL UNIVERSITY
ITHACA, NEW YORK 14853
USA
E-MAIL: rtd1@math.cornell.edu

DEPARTMENT OF MATHEMATICS AND STATISTICS
CONCORDIA UNIVERSITY
MONTREAL, QUEBEC
CANADA H3G 1M8
E-MAIL: lpopovic@mathstat.concordia.ca